\documentclass[12pt]{article}
\usepackage{latexsym,amsmath,amsthm}
\usepackage{amssymb}
\usepackage{graphicx}
\usepackage{appendix}
\usepackage{tikz}
\usepackage[colorlinks=true, linkcolor=blue]{hyperref}
\everymath{\displaystyle}
\usepackage{setspace}

\usepackage{tikz}
\usetikzlibrary{fit,shapes,positioning,calc,decorations.pathreplacing}

\newtheorem{thm}{Theorem}[section]

\newtheorem{cor}[thm]{Corollary}

\newtheorem{claim}[thm]{Claim}

\theoremstyle{definition}

\theoremstyle{definition}
\newtheorem{note}[thm]{Note}

\theoremstyle{definition}
\newtheorem{question}[thm]{Question}

\theoremstyle{definition}

\theoremstyle{definition}

\theoremstyle{definition}

\theoremstyle{definition}

\theoremstyle{definition}
\newtheorem{conv}[thm]{Convention}

\theoremstyle{remark}

\usepackage[margin=1in]{geometry}

\usepackage{url,xcolor}
\usepackage{microtype}

\newcommand{\seqnum}[1]{\href{https://oeis.org/#1}{\rm #1}}

\begin{document}

\title{Independent set sequence of some linear hypertrees} 

\author{David Galvin\thanks{Department of Mathematics, University of Notre Dame, Notre Dame IN, United States; dgalvin1@nd.edu. Research supported in part by the Simons Foundation Gift number 854277.} \and Courtney Sharpe\thanks{Protiviti, Chicago IL, United States; csharpe2@alumni.nd.edu.}}

\maketitle

\begin{abstract}
The independent set sequence of trees has been well studied, with much effort devoted to the (still open) question of Alavi, Malde, Schwenk and Erd\H{o}s on whether the independent set sequence of a tree is always unimodal.

Much less attention has been given to the independent set sequence of hypertrees. Here we study some natural first questions in this realm. We show that the strong independent set sequences of linear hyperpaths and of linear hyperstars are unimodal (actually, log-concave). For uniform linear hyperpaths we obtain explicit expressions for the number of strong independent sets of each possible size, both via generating functions and via combinatorial arguments. We also consider the uniform linear hypercomb with $n$ edges on the spine, and show that its strong independent set sequence is unimodal except possibly for a portion of length $o(n)$.
\end{abstract}

\section{Introduction}

In 1987 Alavi, Malde, Schwenk and Erd\H{o}s considered the {\it vertex independent set sequence}, or simply {\it independent set sequence}, of a graph. This is the sequence $(i_k(G))_{k \geq 0}$ where $i_k(G)$ is the number of vertex independent sets (sets of pairwise non-adjacent vertices) of size $k$ in a graph $G$. (Note that all graphs and hypergraphs in this paper are simple\,---\,all edges have multiplicity one\,---\,and undirected.)

The {\it edge independent set sequence}, better known as the {\it matching sequence}, is the sequence $(m_k(G))_{k \geq 0}$ where $m_k(G)$ is the number of matchings (sets of pairwise disjoint edges) of size $k$ in $G$. A corollary of a seminal result of Heilmann and Lieb \cite{HL1972} is that the matching sequence is always unimodal, that is, that there is a mode $k$ such that
$$
m_0(G) \leq m_1(G) \leq \cdots \leq m_k(G) \geq m_{k+1}(G) \geq \cdots.
$$
The notion of unimodality (and the related ideas of log-concavity and real-rootedness, to be discussed later) are ubiquitous in combinatorics; see for example the surveys \cite{B2015, B1994, S1989}.

Wilf asked the question ``is the independent set sequence of every graph similarly unimodal?''. In \cite{AMSE1987} it was reported that he was skeptical of an affirmative answer, and it turned out that his skepticism was justified. Alavi, Malde, Schwenk and Erd\H{o}s \cite{AMSE1987} showed that the independent set sequence of a graph can fail spectacularly to be unimodal. Specifically they showed that for every $m$ and every permutation $\pi$ of $\{1, \ldots, m\}$ there is a graph $G$ whose largest independent set has size $m$ and for which
$$
i_{\pi(1)}(G) < i_{\pi(2)}(G) < \cdots < i_{\pi(m)}(G).
$$
In other words, the independent set sequence can exhibit arbitrary rises and falls. (Note that we do not consider $i_0(G)$ here as it always takes the smallest possible non-zero value, $1$).

This ``roller coaster'' observation naturally led Alavi, Malde, Schwenk and Erd\H{o}s to consider whether there are infinite families of graphs for which the independent set sequence {\it is} unimodal. Since the matching sequence of a graph is the independent set sequence of its line graph, it follows from Heilmann and Lieb's result that line graphs have unimodal independent set sequences. This was generalized by Hamidoune \cite{H1990}: claw-free graphs (graphs without an induced $K_{1,3}$) have unimodal independent set sequences.      

Alavi, Malde, Schwenk and Erd\H{o}s asked about another very basic family of graphs. It is easy to check that the independent set sequence of a path is unimodal, and even easier to check that the same is true of a star. That unimodality holds for these two most extreme trees provides some reason to expect that it holds for all trees.
\begin{question} 
\label{q-amse}
(\cite{AMSE1987}) Do all trees have unimodal independent set sequence? What about all forests?
\end{question}
The {\it independent set polynomial} $p_G(x)$ is the generating polynomial $\sum_{k \geq 0} i_k(G)x^k$. It is an easy check that if $G$ has components $G_1, \ldots, G_m$ then $p_G(x)=p_{G_1}(x) \cdots p_{G_m}(x)$. Since the product of polynomials with unimodal coefficient sequences need not have a unimodal coefficient sequence, an affirmative answer to Question \ref{q-amse} for trees does not immediately imply an affirmative answer for forests. As an example of this phenomenon of unimodality not behaving well under multiplication, consider the graphs $G_1 = K_{100} + 3K_7$ (the join of a clique on 100 vertices and a disjoint union of three cliques each on 7 vertices) and $G_2=K_{90}+3K_7$. We have
$$
p_{G_1}(x) = 1 + 121x + 147x^2 + 343x^3, ~~~p_{G_2}(x) = 1 + 111x + 147x^2 + 343x^3,
$$
which are both unimodal, but
$$
p_{G_1 \cup G_2}(x) = p_{G_1}(x)p_{G_2}(x) = 1 + 232x + 13725x^2 + 34790x^3 + 101185x^4 + 100842x^5 + 117649x^6,
$$
which is not. (This example is due to Levit and Mandrescu \cite{LM2006}.)
 
Remarkably, Question \ref{q-amse} remains fairly wide open. It has been verified for some infinite, mostly recursively defined, families of trees (see e.g.\ \cite{GH2018} and the references therein). It has also been verified for all trees on 26 or fewer vertices \cite{KLYM2023}. Notably, although all trees on 25 or fewer vertices satisfy the stronger property that their independence polynomials are log-concave (i.e., satisfy $a_k^2 \geq a_{k-1}a_{k+1}$ for all $k \geq 1$), two trees on 26 vertices are found in \cite{KLYM2023} to have unimodal but not log-concave independent set sequences, and in \cite{BR2025, G2025, KL2023, RS2025} infinite families of trees with non-log-concave (though still unimodal) independent set sequences are found. This perhaps explains why Question \ref{q-amse} has proven to be so thorny\,---\,typically it is easier to demonstrate log-concavity than unimodality.

Even answering Question \ref{q-amse} for the uniform random labelled tree (i.e., showing that at least a proportion $1-o(1)$ of all labelled trees on $n$ vertices have unimodal independent set sequences) has proven to be a difficult problem. See \cite{BG2021, H2025} where this question is addressed and only partial results are obtained.

\subsection{Hypertrees}

The goal of this paper is to begin the project of considering Alavi, Malde, Schwenk and Erd\H{o}s' question for {\it hypertrees}. A {\it hypergraph} is a pair $(V,E)$ where $V$ is a set of vertices and $E$ is a set of {\it edges}, which are subsets of $V$. If all subsets have size $2$ then a hypergraph is simply a graph. A hypergraph $H$ is a {\it hypertree} if there is a tree ${\rm tree}(H)$ on vertex set $V(H)$ with the property that for each $e \in E(H)$ the subgraph of ${\rm tree}(H)$ induced by $e$ is connected. Every tree $T$ is a hypertree\,---\,just take ${\rm tree}(T)=T$. 

Before asking any question about independent sets in hypertrees, we need to clarify what we mean, since as is almost always the case when moving from graphs to hypergraphs there are competing definitions. A {\it weak independent set} in a hypergraph is a subset of vertices that does not contain all the vertices of any edge of the hypergraph, while a {\it strong independent set} is one that contains at most one vertex from each edge. Note that these two notions coincide for a graph.
\begin{conv}
In this paper we work exclusively with strong independent sets, and will simply refer to them as {\it independent sets} from here on. 
\end{conv}

It is not the case that every hypertree has unimodal independent set sequence. Indeed, consider the hypergraph on vertex set $\{v_1, \ldots, v_{10}, w_1, \ldots, w_6\}$ with $\{v_1,\ldots, v_{10}, w_i\}$ an edge for each $i$. Independent sets of size $2$ or greater in this hypergraph  cannot include any vertices from among $\{v_1, \ldots, v_{10}\}$ but can include vertices from $\{w_1, \ldots, w_6\}$ without restriction. It follows that the independent set sequence of the hypergraph is $(1,16,15,20,15,6,1)$, which is not unimodal.   

\begin{figure}[htp]
    \centering
    \includegraphics[width=11cm]{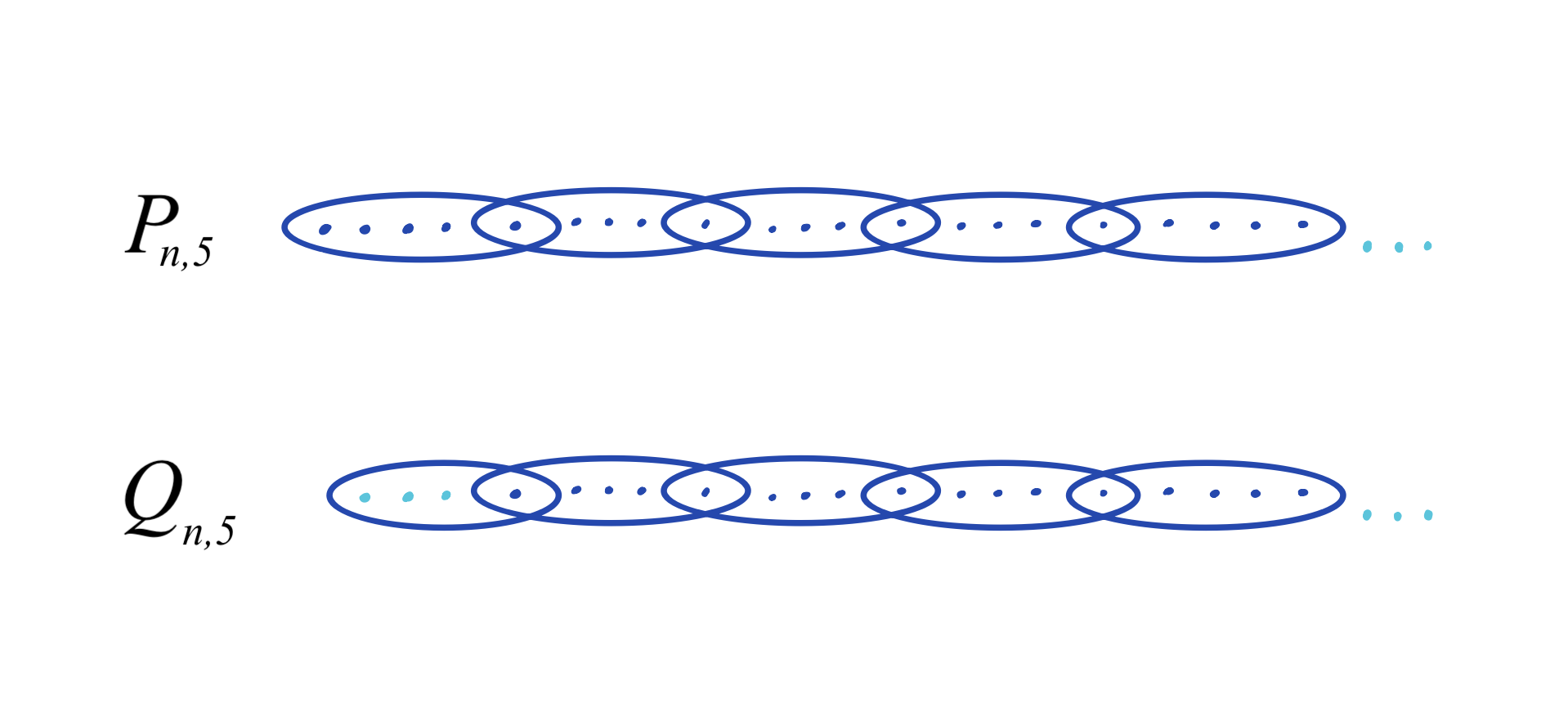}
    \caption{The first 5 edges of the linear hyperpaths $P_{n,5}$ (top) and $Q_{n,5}$ (bottom). Note that the first edge of $Q_{n,5}$ has one vertex fewer than the first edge of $P_{n,5}$.}
    \label{fig:Q1}
\end{figure}

To remove from consideration examples like this, our focus in this paper will be on linear hypertrees. A hypergraph is {\it linear} if any two edges intersect in at most one vertex. We begin by considering two infinite families of linear hypertrees, namely the linear hyperpaths and the linear hyperstars. Let integers $n \geq 1$ and $s_1, s_2, \ldots, s_n \geq 2$ be given.
\begin{itemize}
\item The {\it linear hyperpath} $P(s_1, \ldots, s_n)$ is the hypergraph with $n$ edges $e_1, \ldots, e_n$, with $|e_i|=s_i$ for each $i$, with $e_i$ and $e_{i+1}$ sharing a single vertex for each $i<n$, and with no other pairs of edges sharing a vertex. 
\item If all $s_i=\ell$ for some $\ell \geq 2$ then we use the notation $P_{n,\ell}$, and refer to this hypergraph as the {\it $\ell$-uniform, $n$-edge, linear hyperpath} (see Figure \ref{fig:Q1}, which also shows an auxiliary hypergraph that will be defined later).
\item The {\it linear hyperstar} $S(s_1, \ldots, s_n)$ is the hypergraph with $n$ edges $e_1, \ldots, e_n$, with $|e_i|=s_i$ for each $i$, and with a single vertex (the {\it center} of the star) in common to all of the edges.  
\end{itemize}

Note that $P_{n,2}$ is just the (graph) path on $n$-edges, while $S(2, 2,\ldots, 2)$ is the (graph) star on $n$-edges. Note also that we use slightly non-standard notation here, with $P_{n,2}$ being the path on $n$ edges (and so $n+1$ vertices) rather than on $n$ vertices (and so $n-1$ edges). We adopt this convention because it is better suited to the uniform hypergraph ($\ell \geq 3$) setting.

Just as with trees, the unimodality of independent set sequences of linear hyperpaths and linear hyperstars is fairly straightforward. Before stating the result, we introduce some further concepts. A sequence $(a_k)_{k=0}^n$ of real numbers is {\it log-concave} if for all $1 \leq k \leq n-1$ it holds that $a_k^2 \geq a_{k-1}a_{k+1}$. If a log-concave sequence is non-negative and has no internal zeroes (meaning that there is no triple $0 \leq i \leq j \leq k \leq n$ with $a_i >0$, $a_j=0$ and $a_k >0$) then it is easy to check that the sequence is unimodal. Often log-concavity of a sequence is easier to establish than unimodality, since to verify unimodality one must identify the location of the mode, but there is no such difficulty when establishing log-concavity.

There is a stronger property of finite real sequences than log-concavity. A sequence $(a_k)_{k=0}^n$ is {\it real-rooted} if the polynomial $\sum_{k=0}^n a_kx^k$ has only real roots. A  result that ultimately goes back to Newton is that if a non-negative sequence is real-rooted, then it is log-concave, and in fact satisfies the stronger inequalities
\begin{equation} \label{eq-ulc}
a_k^2 \geq \left(1+\frac{1}{k}\right)\left(1+\frac{1}{n-k}\right)a_{k-1}a_{k+1}
\end{equation}
for all $1 \leq k \leq n-1$ (see e.g. \cite[Section 9.3]{B2016}). We note incidentally that for a non-negative sequence with no internal zeros neither of the implications in the chain  
\begin{center}
real-rooted $\Rightarrow$ log-concave $\Rightarrow$ unimodal
\end{center}
can be reversed. 

We now return to hyperstars and hyperpaths. Since the independent set sequences of $P(s_1, \ldots, s_n)$ and $S(s_1, \ldots, s_n)$ have no internal zeroes, our first result\,---\,Theorem \ref{thm-hyperstar,hyperpath} below\,---\,implies that these sequences are unimodal. The proof of Theorem \ref{thm-hyperstar,hyperpath} is given in Section \ref{sec-nonuniform-proofs}.
\begin{thm} \label{thm-hyperstar,hyperpath}
For $n \geq 1$ and $s_1, s_2, \ldots, s_n \geq 2$
\begin{itemize}
\item the independent set sequence of $P(s_1, \ldots, s_n)$ is real-rooted, and
\item the independent set sequence of $S(s_1, \ldots, s_n)$ is log-concave.
\end{itemize}
\end{thm}
 
\begin{note} \label{note-CS}
In a draft of this paper \cite[version 1]{GS2024-v1} we deduced the real-rootedness of the independent set sequence of $P_{n,\ell}$ using an interlacing argument. After we posted the draft on the arXiv preprint server, Ferenc Bencs pointed out to us that we could use results of Chudnovsky and Seymour \cite{CS2007} and Hamidoune \cite{H1990} obtain the more general result presented here.  
\end{note}

Just as the unimodality of the independent set sequences of paths and stars (the trees with the largest and smallest diameters, respectively) lends credence to Alavi, Malde, Schwenk and Erd\H{o}s' speculation that all trees have unimodal independent set sequence, Theorem \ref{thm-hyperstar,hyperpath} makes it plausible to speculate that maybe the (strong) independent set sequence of every linear hypertree is unimodal. We will return to this briefly at the end of the introduction.     

\medskip

The bulk of this paper is concerned with $P_{n,\ell}$, the linear {\it uniform} hyperpath. While the unimodality of the independent set sequence of $P_{n,\ell}$ follows from Theorem \ref{thm-hyperstar,hyperpath}, there is much interesting combinatorics to explore related to explicit formulae and recurrence relations for the terms of the sequence. For $\ell=2$ (graph paths) all is quite straightforward. Let $p_{n,2}^k$ be the number of independent sets of size $k$ in $P_{n,2}$. It is easy to generate a recurrence relation that expresses $p_{n,2}^k$ in terms of $p_{n',2}^{k'}$'s with $n'+k'<n+k$. It is also easy to solve these recurrences (both inductively and combinatorially) to show that (modulo some initial conditions) 
\begin{equation} \label{path-count}
p_{n,2}^k=\binom{n-k+2}{k}.
\end{equation}
Things get more involved when $\ell \geq 3$. Let $p_{n,\ell}^k$ be the number of independent sets of size $k$ in $P_{n,\ell}$. As we will see in the sequel, the natural recurrence relation to calculate $p_{n,\ell}^k$ requires introducing an auxiliary (non-uniform) hyperpath, and is in fact a pair of coupled recurrence relations. From these recurrences it is possible to conjecture an explicit but non-obvious formula for $p_{n,\ell}^k$, which turns out to be a sum of products of binomial coefficients. It is fairly straightforward (if a little messy) to prove this explicit formula by induction, but it is much less straightforward to give a direct combinatorial argument. Finally, it is easy to decouple the pair of recurrence relations to produce an expression for $p_{n,\ell}^k$ in terms of $p_{n',\ell}^{k'}$'s with $n'+k'<n+k$, but it is by no means obvious how to directly (combinatorially) verify this expression.

In Sections \ref{subsec-induction-path}, \ref{subsec-comb-path} and \ref{subsec-fib-rec} we give the proof of the following result, that makes precise what was alluded to in the previous paragraph. 
\begin{thm} \label{thm-pknl-exact}
For all $\ell \geq 2$, we have
\begin{itemize}
\item $p^k_{0,\ell}=1$ if $k=0$ and $p^k_{0,\ell}=0$ if $k>0$,
\item $p^0_{1,\ell}=1$, $p^1_{1,\ell}=\ell$ and $p^k_{1,\ell}=0$ for $k \geq 2$,
\item $p^0_{n,\ell}=1$ for $n \geq 0$,
\item $p^1_{n,\ell}=0$ for $n = 0$ and $p^1_{n,\ell}=n\ell-(n-1)$ for $n > 0$, and
\item for $n \geq 2$ and $k \geq 2$ 
\begin{equation} \label{pknl-summation-formula}
p^k_{n,\ell}=(\ell-1)^2\sum_{j=0}^{k-2}  (\ell - 2)^{k-j-2} \binom{k-2}{j}\binom{n-j}{k}.
\end{equation}
\end{itemize}
We also have, for $n \geq 3$, $k \geq 1$ and $\ell \geq 2$, the recurrence  
\begin{equation} \label{fib-rec}
p^k_{n,\ell} = p^k_{n-1,\ell} + p^{k-1}_{n-2,\ell} + (\ell-2)p^{k-1}_{n-1,\ell}.
\end{equation}
\end{thm}
Notice that when $\ell=2$ \eqref{pknl-summation-formula} reduces to \eqref{path-count} (as long as we interpret $0^0=1$), and that \eqref{fib-rec} reduces to the familiar (and easy to derive) recurrence $p^k_{n,2} = p^k_{n-1,2} + p^{k-1}_{n-2,2}$. 

In Section \ref{subsec-induction-path} we begin with the easy algebraic derivation of \eqref{fib-rec}, and then describe how \eqref{pknl-summation-formula} may be conjectured heuristically. Once the correct formula has been obtained, it is fairly straightforward (if a little messy) to prove it by induction. In Section \ref{subsec-comb-path} we give a more satisfying combinatorial explanation for \eqref{pknl-summation-formula}. In Section \ref{subsec-fib-rec} we give a combinatorial explanation for \eqref{fib-rec}.

The algebraic derivation of \eqref{fib-rec} will be through the independence polynomial of $P_{n,\ell}$, that is, through 
$$
P_{n,\ell}(x) = \sum_{k=0}^n p^k_{n,\ell} x^k.
$$
(Note that $P_{n,\ell}$ has $n$ edges, so $p^k_{n,\ell}=0$ for $k > n$.) We will establish the recurrence
\begin{equation} \label{eq-Px-pure}
P_{n,\ell}(x) = (1+(\ell-2)x)P_{n-1,\ell}(x) + xP_{n-2,\ell}(x)
\end{equation}
for $n \geq 3$ (with initial conditions $P_{0,\ell}(x)=1$, $P_{1,\ell}(x)=1+\ell x$ and $P_{2,\ell}(x)=1+(2\ell-1)x+(\ell-1)^2x^2$) from which \eqref{fib-rec} follows immediately by equating coefficients of $x^k$ on both sides. Observe that we could have chosen instead to declare $P_{0,\ell}(x)=1+x$. This would extend the validity of \eqref{eq-Px-pure} to $n \geq 2$, but at the cost of introducing a very unnatural initial condition (there is no $\ell$-uniform hypergraph for $\ell \geq 3$ with independence polynomial $1+x$).    

\begin{note} \label{note-comb}
In a draft of this paper \cite[version 1]{GS2024-v1} we derived \eqref{fib-rec} algebraically but left a combinatorial explanation as a problem, since we had been unable to find one. After we posted the draft on the arXiv preprint server Ferenc Bencs communicated to us the combinatorial proof that we present in Section \ref{subsec-fib-rec}.  
\end{note}

\begin{figure}[htp]
    \centering
    \includegraphics[width=11cm]{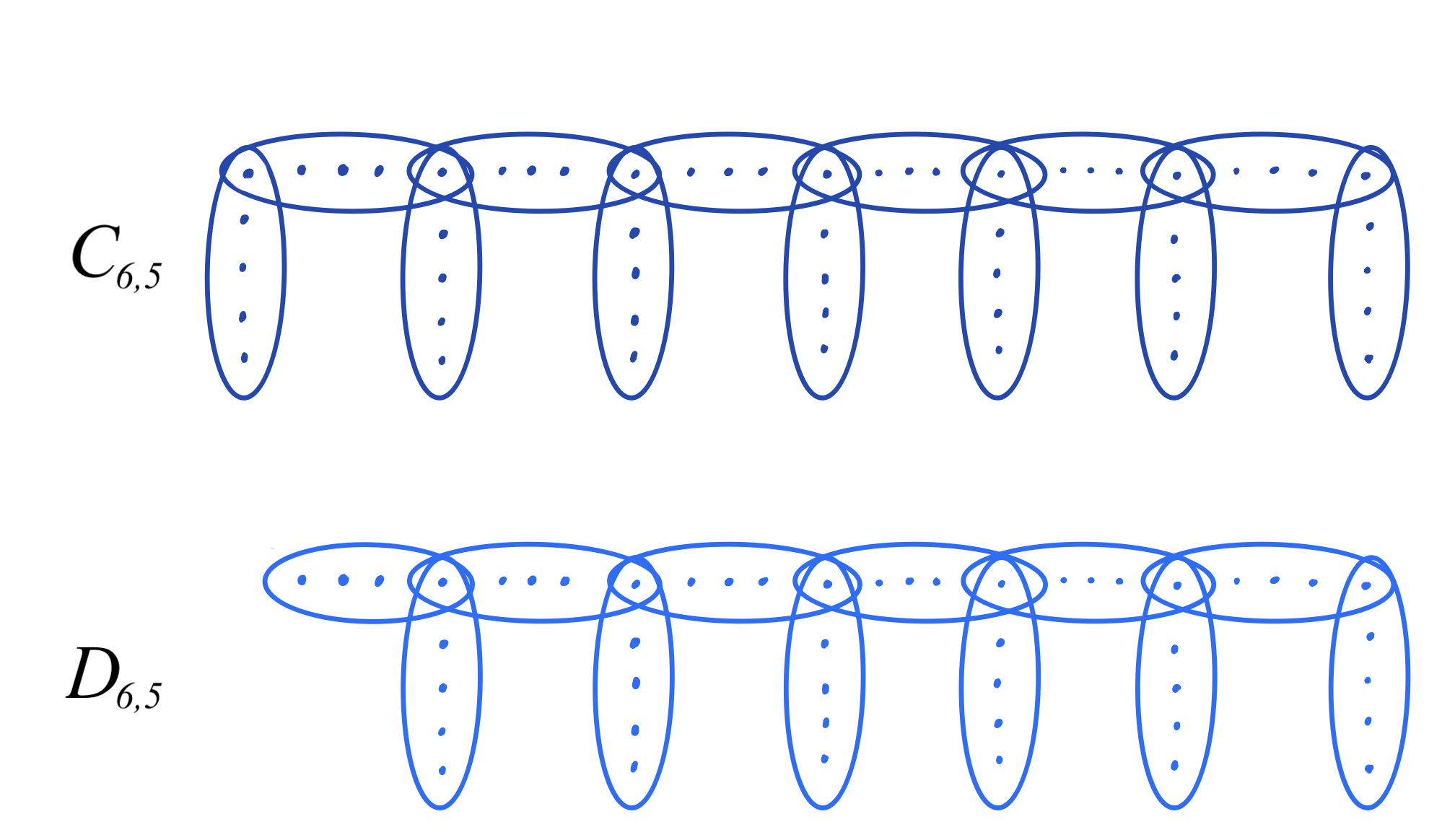}
    \caption{The linear uniform hypercomb $C_{6,5}$ (top) and the auxiliary hypercomb $D_{6,5}$ (bottom).}
\label{fig:CandD}
\end{figure}

Towards the goal of making plausible the speculation that linear hypertrees have unimodal independent set sequences, we also consider another infinite family of hypertrees, one for which the results of Chudnovsky and Seymour and of Hamidoune alluded to in Note \ref{note-CS} cannot be applied. This is the family of uniform linear hypercombs. For $n \geq 1$ the {\it uniform linear hypercomb} $C_{n,\ell}$ is the hypergraph obtained from $P_{n,\ell}$ by adding, 
\begin{itemize}
\item for each vertex $v$ that is in two edges of $P_{n,\ell}$, an edge of size $\ell$ that includes $v$ and $\ell-1$ new vertices, and also 
\item for one vertex $w$ in the first edge of $P_{n,\ell}$ but not in any other edges, and for one vertex $w'$ in the last edge of $P_{n,\ell}$ but not in any other edges, adding an edge of size $\ell$ that includes $w$ and $\ell-1$ new vertices, and an edge of size $\ell$ that includes $w'$ and $\ell-1$ new vertices. If $n\geq 2$ then necessarily $w \neq w'$; if $n=1$ then we enforce that $w \neq w'$.
\end{itemize}
See Figure \ref{fig:CandD}, which also shows an auxiliary hypergraph that will be defined later. 
When $\ell=2$ this graph is referred to as a {\em centipede} or a {\em comb}. That the centipede has unimodal independent set sequence was first established by Levit and Mandrescu in \cite{LM2002}, and later Zhu \cite{Z2007} showed that the sequence is real-rooted. 

We have observed computationally that the independent set sequence of $C_{n,\ell}$ is unimodal for all $n \leq 500$ and $\ell \leq 12$. Here we establish that for all $n$ and $\ell$ the sequence is unimodal expect (possibly) for an $O(\sqrt{n})$ portion. In Theorem \ref{thm-comb} and Corollary \ref{cor-comb} below we denote by $c^k_{n,\ell}$ the number of independent sets in $C_{n,\ell}$ of size $k$.

\begin{thm} \label{thm-comb}
For $\ell \geq 3$, the independence polynomial $C_{n,\ell}(x)$ of $C_{n,\ell}$ satisfies the recurrence
\begin{itemize}
\item $C_{0,\ell}(x) = 1+\ell x$ (an initial condition),
\item $C_{1,\ell}(x) = 1+(3\ell-2)x+(2(\ell-1)(\ell-2)+(\ell-1)^2+2(\ell-1))x^2+(\ell-1)^2(\ell-2)x^3$, and
\item for $n \geq 2$
\begin{equation} \label{eq-comb-rec}
C_{n,\ell}(x) = (1+(\ell-1)x)\left((1+(\ell-2)x)C_{n-1,\ell}(x)+xC_{n-2,\ell}(x)\right).
\end{equation}
\end{itemize}
\end{thm}
The proof of Theorem \ref{thm-comb}, and of Corollary \ref{cor-comb} below, are given in Section \ref{subsec-comb-proofs}.  
\begin{cor} \label{cor-comb}
For each $\ell \geq 3$ there are constants $K_\ell, L_\ell >0$ such that the independent set sequence $(c^k_{n,\ell})_{k=0}^{2n+1}$ of $C_{n,\ell}$ is weakly increasing up to $K_\ell n - L_\ell\sqrt{n}$ (i.e., $c^{k-1}_{n,\ell} \leq c^k_{n,\ell}$ for $k \leq K_\ell n - L_\ell\sqrt{n}$) and is weakly decreasing from $K_\ell n + L_\ell\sqrt{n}$ (i.e., $c^k_{n,\ell} \geq c^{k+1}_{n,\ell}$ for $k \geq K_\ell n + L_\ell\sqrt{n}$). 
\end{cor}

\begin{note} \label{note-hcomb}
In a draft of this paper \cite[version 3]{GS2024-v1} we claimed that the independent set sequence of $C_{n,\ell}$ is unimodal (in fact, log concave) for all $n$ and $\ell$. During the refereeing process we discovered an error in our proof, which has led us to the weaker statement. Our approach to Corollary \ref{cor-comb} was partly inspired by ideas presented in \cite{WZ2011}, one of the references in \cite{BRGGGS2019}; we thank Ferenc Bencs for bringing \cite{BRGGGS2019} to our attention. 
\end{note}

In proving Theorem \ref{thm-comb} and Corollary \ref{cor-comb} we work in a more general setting, attaching a copy of the same fixed hypergraph $G$ to the first vertex of $P_{n,\ell}$, the last vertex of $P_{n,\ell}$, and all vertices that are in two edges of $P_{n,\ell}$ (for $C_{n,\ell}$, $G$ is a single edge of order $\ell$). Our approach to studying the independent set sequence of this family of hypergraphs is very much in the spirit of \cite[Theorem 3.1 and subsequent applications]{WZ2011}. We defer the details to Section \ref{subsec-comb-proofs}.

\medskip

We end with some questions for further consideration.
\begin{question}
\textcolor{white}{X}
\begin{enumerate}
\item Is the strong independent set sequence of a linear hypertree unimodal?
\item If not, then what about a linear uniform hypertree?
\item What if ``strong independent set'' is replaced by ``weak independent set''?
\item Can the gap of length $O(\sqrt{n})$ in Corollary \ref{cor-comb} be shortened or eliminated? 
\end{enumerate}
\end{question}

Much of the progress to date on Alavi, Malde, Schwenk and Erd\H{o}s' tree question has consisted of establishing unimodality of the independent set sequence of trees that lie in particular (often recursively defined) families. One might tackle analogous families in the hypertree setting, although the fact that (as seen in this paper) some complexity arises even in the simple families of linear uniform hyperpaths and linear uniform hypercombs suggests that the results in the hypertree setting maybe be harder than their counterparts in the tree setting.

\medskip

Returning to one of the main focuses of this note, the linear uniform hyperpath: our combinatorial derivation of the identity \eqref{pknl-summation-formula} does not make it obvious that $p^k_{n,\ell}$ is a multiple of $(\ell-1)^2$ when $n, k, \ell \geq 2$. Is there a simple explanation for this phenomenon? 

\section{Proofs} \label{sec-proofs}

\subsection{Proof of Theorem \ref{thm-hyperstar,hyperpath}} \label{sec-nonuniform-proofs}

We begin with the hyperpath $P(s_1, \ldots, s_n)$. Let $GP(s_1, \ldots, s_n)$ be the graph obtained from $P(s_1, \ldots, s_n)$ by replacing each edge $e_i$ with a complete graph on the same set of vertices. Independent sets in $GP(s_1, \ldots, s_n)$ are exactly (strong) independent sets in $P(s_1, \ldots, s_n)$. Since $GP(s_1, \ldots, s_n)$ is claw-free (has no induced copy of $K_{1,3}$), it follows from Hamidoune's result \cite{H1990} cited earlier that the independent set sequence of $P(s_1, \ldots, s_n)$ is log-concave. Hamidoune's result was extended by Chudnovsky and Seymour \cite{CS2007} who showed that the independent set sequence of a claw-free graph is real-rooted, and so in fact the independent set sequence of $P(s_1, \ldots, s_n)$ is real-rooted.   

We now move on to the linear hyperstar. Noting that other than the independent set of size $1$ that consists of the center of the star (vertex $v$, say), an independent set in $S(s_1,s_2, \ldots, s_n)$ is the union of arbitrary subsets of size $0$ or $1$ of $e_i\setminus\{v\}$ as $i$ runs between $1$ and $n$, we see that the independence polynomial of  $S(s_1,s_2, \ldots, s_n)$ is 
$$
x + \prod_{i=1}^n (1+(s_i-1)x).
$$
 
For $n=1, 2$ this polynomial is easy seen to be real-rooted and so has log-concave coefficient sequence. For $n\geq 3$ the independence polynomial of $S(s_1, \ldots, s_n)$ is not always real-rooted (consider $n=3$ and $s_1=s_2=s_3=2$, for example), but we can still establish log-concavity of the coefficient sequence.   

Let $a_0, \ldots, a_n$ be defined by
$$
\prod_{i=1}^n (1+(s_i-1)x) = \sum_{j=0}^n a_jx^j,
$$
so that the independent set sequence of $S(s_1, \ldots, s_n)$ is $(a_0,1+a_1,a_2, \ldots, a_n)$. Because $\prod_{i=1}^n(1+(s_i-1)x)$ is real-rooted its coefficient sequence is log-concave, and this automatically gives us all of the log-concavity relations for the independent set sequence of $S(s_1, \ldots, s_n)$, except
$$
a_2^2 \geq (a_1+1)a_3~~~\mbox{and}~~~(a_1+1)^2 \geq a_0a_3.
$$
The second of these follows immediately from $a_1^2 \geq a_0a_3$. For the first we use that in fact the coefficient sequence of $\prod_{i=1}^n(1+(s_i-1)x)$ satisfies the Newton inequalities \eqref{eq-ulc}, so that
$$
a_2^2 \geq \left(\frac{3n-3}{2n-4}\right)a_1a_3.
$$
The relation $a_2^2 \geq (a_1+1)a_3$ would thus follow from
$$
\left(\frac{3n-3}{2n-4}\right)a_1a_3 \geq (a_1+1)a_3
$$
or equivalently $a_1 \geq (2n-4)/(n+1)$. Since all $s_i \geq 2$ and $n \geq 3$ we have $a_1\geq 3$, and since $3 \geq (2n-4)/(n+1)$ for all $n \geq 3$ we indeed get $a_1 \geq (2n-4)/(n+1)$. We conclude that the independent set sequence of $S(s_1, \ldots, s_n)$ is log-concave.  

\subsection{Deriving a formula for $p_{n,\ell}^k$} \label{subsec-induction-path}

For ordinary paths there is a very simple recurrence for calculating $p^k_{n,2}$ (once suitable initial conditions have been established), namely
$$
p^k_{n,2} = p^k_{n-1,2} + p^{k-1}_{n-2,2}, 
$$
obtained by considering whether or not the first vertex of the path is in the independent set. For $\ell \geq 3$ we need to introduce an auxiliary hypergraph, which we denote by $Q_{n,\ell}$, that is obtained from the hyperpath $P_{n,\ell}$ by removing the vertex $v_1$ (so $Q_{n,\ell}$ has one edge of size $\ell-1$, with the rest having size $\ell$; see Figure \ref{fig:Q1}). Let the edges of $Q_{n,\ell}$ be $e_1', e_2, \ldots, e_n$, and denote by $q^k_{n,\ell}$ the number of independent sets of size $k$ in $Q_{n,\ell}$. 

There are a pair of coupled recurrences for $p^k_{n,\ell}$ and $q^k_{n,\ell}$, obtained by considering whether or not the independent set includes a vertex from $e_1\setminus \{v_{\ell}\}$ or $e'_1\setminus \{v_{\ell}\}$. Indeed, for $n, k \geq 2$ we have 
\begin{equation} \label{recurrence}
\begin{array}{rcl}
p^k_{n,\ell} & = & p^k_{n-1,\ell} + (\ell-1)q^{k-1}_{n-1,\ell} \\
& \mbox{and} & \\
q^k_{n,\ell} & = & p^k_{n-1,\ell} + (\ell-2)q^{k-1}_{n-1,\ell}.
\end{array}
\end{equation}
Observe that when $\ell=2$ the second relation above reduces to $q^k_{n,\ell} = p^k_{n-1,\ell}$, and the system above is easily seen to reduce to the standard recurrence for $p_{n,2}^k$.

The initial conditions for the recurrences for $p^k_{n,\ell}$ and $q^k_{n,\ell}$ are:
\begin{itemize}
\item For $n=0$: $p^0_{0,\ell}=q^0_{0,\ell}=1$ (due to the empty set) and $p^k_{0,\ell}=q^k_{0,\ell}=0$ for $k \geq 1$.
\item For $n=1$: $p^0_{1,\ell}=q^0_{1,\ell}=1$, $p^1_{1,\ell}=\ell$, $q^1_{1,\ell}=\ell-1$, and $p^k_{1,\ell}=q^k_{1,\ell}=0$ for $k \geq 2$.
\item For $k=0$: $p^0_{n,\ell}=q^0_{n,\ell}=1$.
\item For $k=1$: $p^1_{n,\ell}=n\ell-(n-1)$ and $q^1_{n,\ell}=n\ell-n$ for $n \geq 2$. 
\end{itemize}  

The recurrences \eqref{recurrence} for $p^k_{n,\ell}$ and $q^k_{n,\ell}$ can be decoupled to allow $p^k_{n,\ell}$ to be expressed in terms of $p^{k'}_{n',\ell}$ with $k'+n' < k+n$ and with no reference to $q^{k'}_{n',\ell}$'s. To achieve this decoupling it is easiest to work with the independence polynomials of $P_{n,\ell}$ and $Q_{n,\ell}$. Setting
$$
P_{n,\ell}(x) = \sum_{k \geq 0} p^k_{n,\ell}x^k ~\mbox{and}~Q_{n,\ell}(x) = \sum_{k \geq 0} q^k_{n,\ell}x^k
$$ 
we have (as before by first considering those independent sets in which none of the vertices unique to the first edge are in the independent set, and then considering those in which one such vertex is in the independent set)
\begin{equation} \label{GPrecurrence}
P_{n,\ell}(x)=P_{n-1,\ell}(x) + (\ell-1)xQ_{n-1,\ell}(x)
\end{equation}
and
\begin{equation} \label{GQrecurrence}
Q_{n,\ell}(x)=P_{n-1,\ell}(x) + (\ell-2)xQ_{n-1,\ell}(x)
\end{equation}
for $n \geq 2$, with initial conditions 
$$
P_{0,\ell}(x)=Q_{0,\ell}(x)=1, P_{1,\ell}(x)=1+\ell x~\mbox{and}~
Q_{1,\ell}(x)=1+(\ell-1) x.
$$
From \eqref{GPrecurrence} we have
$$
Q_{n-1,\ell}(x) = \frac{P_{n,\ell}(x)-P_{n-1,\ell}(x)}{(\ell-1)x}~~~\mbox{and}~~~Q_{n,\ell}(x) = \frac{P_{n+1,\ell}(x)-P_{n,\ell}(x)}{(\ell-1)x}.
$$
Inserting these into  \eqref{GQrecurrence} yields
$$
\frac{P_{n+1,\ell}(x)-P_{n,\ell}(x)}{(\ell-1)x} = P_{n-1,\ell}(x) + \frac{(\ell-2)x\left(P_{n,\ell}(x)-P_{n-1,\ell}(x)\right)}{(\ell-1)x}.
$$
After rearranging terms and shifting indices we get that for $n \geq 3$,
$$
P_{n,\ell}(x) = (1+(\ell-2)x)P_{n-1,\ell}(x) + xP_{n-2,\ell}(x).
$$
Combining this with the initial conditions $P_{0,\ell}(x)= 1$ and $P_{1,\ell}(x)= 1+\ell x$ (as introduced above) and
$$
P_{2,\ell}(x) = 1 + (2\ell -1)x + (\ell-1)^2 x^2
$$
(an easy direct count) we have completely decoupled $P_{n,\ell}(x)$ from $Q_{n,\ell}(x)$. Extracting the coefficient of $x^n$ from both sides of \eqref{eq-Px-pure} we obtain \eqref{fib-rec}. 

To conjecture an explicit formula for $p_{n,\ell}^k$ for $n, k \geq 2$, we generate values using the recurrence relation, and consult the \href{https://oeis.org}{On-Line Encyclopedia of Integer Sequences} ({\it OEIS}) \cite{OEIS}. For $\ell=3$, the initial part of the array $(p_{n,\ell}^k)_{n, k \geq 0}$ is shown in Table \ref{table-pkn3} (blank entries are $0$).

\begin{table}[ht!]
\begin{center}
\begin{tabular}{r|rrrrrrrc}
$p^k_{n,3}$ & $k=0$ & 1 & 2 & 3 & 4 & 5 & 6 & $\cdots$\\
\hline
$n=0$ & 1 \\
1 &     1 & 3 \\
2 &     1 & 5   & 4   &  \\
3 &     1 & 7   & 12  & 4  & \\
4 &     1 & 9  & 24  & 20  & 4 \\
5 &     1 & 11  & 40  & 56  & 28   & 4 \\
6 &     1 & 13 & 60 & 120 & 104 & 36  & 4 \\
$\vdots$ & $\vdots$ & $\vdots$ & $\vdots$ & $\vdots$ & $\vdots$ & $\vdots$ & $\vdots$ & $\ddots$ 
\end{tabular}
\caption{The count of independent sets of size $k$ in $P^k_{n,3}$.} \label{table-pkn3}
\end{center}
\end{table}

All entries from the $k=2$ column on appear to be multiples of $4$. Table \ref{table-gf} shows the most likely contender from OEIS for the $k$th column of the array in Table \ref{table-pkn3} ($k \geq 2$) once this factor of 4 is pulled out. Table \ref{table-gf} also shows each sequence's ordinary generating function (as given at OEIS).     

\begin{table}[ht!]
\begin{center}
\begin{tabular}{ | c | c | c |  c |}
\hline
$k$ & Sequence & Likely OEIS match & Generating Function \\ 
\hline
$2$ & $(1, 3, 6, 10, ...)$ & \seqnum{A000217}  &  $\frac{x}{(1-x)^3}$ \\ \hline
$3$ & $(1, 5, 14, 30, ...)$ & \seqnum{A000330} & $\frac{x(1+x)}{(1-x)^4}$ \\ \hline
$4$  & $(1, 7, 26, 70, ...)$ 
 & \seqnum{A006325} & $\frac{-x^2(x+1)^2}{(x-1)^5}$ \\ \hline
$5$ & $(1, 9, 42, 138, ...)$ & \seqnum{A061927}  & $\frac{x(1+x)^3}{(-1+x)^6}$ \\ \hline
$\vdots$ & $\vdots$ & $\vdots$ & $\vdots$ \\
\end{tabular}
\caption{The generating functions of each column sequence from Table \ref{table-pkn3}, for $P_{n,3}$.} 
\label{table-gf}
\end{center}
\end{table}

Adjusting each generating function by multiplying by an appropriate power of $x$ (to make sure that the non-zero terms of the associated sequence begin at the correct place), we are led to conjecture that for $k \geq 2$ we have
$$
\sum_{n \geq k} p^k_{n,3}x^n = \frac{4x^k(1+x)^{k-2}}{(1-x)^{k+1}}.
$$
This process can be repeated for larger values of $\ell$; the results for $\ell \leq 6$ are shown in Table \ref{table-lconjectures}. 

\begin{table}[ht!]
\begin{center}
\begin{tabular}{ | c | c | c |}
\hline
$\ell$ & Conjectured Generating Function \\ 
\hline
$3$ & $\sum_{n \geq k} p^k_{n,3}x^n = \frac{4x^k(1+x)^{k-2}}{(1-x)^{k+1}}$ \\ \hline
$4$ & $\sum_{n \geq k} p^k_{n,4}x^n = \frac{9x^k(2+x)^{k-2}}{(1-x)^{k+1}}$ \\ \hline
$5$ & $\sum_{n \geq k} p^k_{n,5}x^n = \frac{25x^k(3+x)^{k-2}}{(1-x)^{k+1}}$ \\ \hline
$6$ & $\sum_{n \geq k} p^k_{n,6}x^n = \frac{36x^k(4+x)^{k-2}}{(1-x)^{k+1}}$ \\ \hline
$\vdots$ & $\vdots$ \\
\end{tabular}
\caption{The conjectured generating functions for each $P_{n,\ell}$, found in patterns from $\ell$-triangle arrays.} 
\label{table-lconjectures}
\end{center}
\end{table}

All of these data motivate the conjecture that
for each $k \geq 2$ and $\ell \geq 3$, we have
\begin{equation} \label{p-gf}
\sum_{n \geq k} p^k_{n,\ell}x^n = \frac{(\ell-1)^2x^k((\ell -2)+x)^{k-2}}{(1-x)^{k+1}}.
\end{equation}

Using the generating function identity
$$
\frac{1}{(1-x)^{k+1}} = \binom{k}{k} + \binom{k+1}{k}x + \binom{k+2}{k}x^2 + \binom{k+3}{k}x^3 +  \cdots
$$
for the denominator above, and the binomial theorem for the numerator, we can extract the coefficient of $x^n$ from the right-hand side of \eqref{p-gf}, heuristically leading us to  \eqref{pknl-summation-formula}.

We now proceed to verify the correctness of \eqref{pknl-summation-formula} by  induction on $n+k$, utilizing the recurrence \eqref{fib-rec}. Since we are claiming validity of \eqref{pknl-summation-formula} for $n, k \geq 2$, there is nothing to do for $k\leq 1$ (and any $n \geq 0$) or $n \leq 1$ (and any $k \geq 0$). It will turn out to be useful to dispense with the case $k=2$ (and $n \geq 2$), as well as $n=2, 3$ (and $k \geq 3$), before beginning the induction on $n+k$.    

First we consider $k=2$. For an independent set of size $2$ drawn from $P_{n,\ell}$ let $i$ be the smallest index such that $e_i$ contains a vertex of the independent set, and let $v$ be that vertex. If $v \in e_1$ but $v$ is not the unique vertex of $e_1$ that is also in $e_2$ ($\ell-1$ options for $v$) then there are $(\ell-1)(n-1)$ options for the second vertex of the independent set. If $v \in e_2\setminus e_3$ (again $\ell-1$ options) then there are $(\ell-1)(n-2)$ options for the second vertex. Continuing the count in this manner we see that
$$
p^2_{n,\ell} = (\ell-1)\sum_{j=1}^{n-1} (\ell-1)(n-j) = (\ell-1)^2 \binom{n}{2}, 
$$
as predicted by \eqref{pknl-summation-formula}.

Next we consider $n=2$. For $k \geq 3$ we have that $p^k_{2,\ell}=0$ (since a $2$-edge hypergraph cannot support an independent set of size greater than $2$), and this is exactly what \eqref{pknl-summation-formula} predicts, owing to the factor of $\binom{2-j}{3}$ (that evaluates to $0$) in every summand.

Similarly for $n=3$ and $k \geq 4$, \eqref{pknl-summation-formula} gives the correct value (namely $0$) for $p^k_{3,\ell}$, while it is easy to see that $p^3_{3,\ell}=(\ell-1)^2(\ell-2)$, as predicted by \eqref{pknl-summation-formula}. 

We can now begin the induction argument. All pairs $(n,k)$ with $n \leq 3$ and $k \leq 2$ have already been considered, so we assume $n+k\geq 7$, $n \geq 4$ and $k \geq 3$. In this range not only is \eqref{fib-rec} valid, but also we may use \eqref{pknl-summation-formula} on all terms on the right-hand side of \eqref{fib-rec} (recall that we are inducting on $n+k$). So we have
\begin{eqnarray}
p^k_{n,\ell} & = & p^k_{n-1,\ell} + p^{k-1}_{n-2,\ell} + (\ell-2)p^{k-1}_{n-1,\ell} \nonumber \\
& = & (\ell-1)^2\sum_j  (\ell - 2)^{k-j-2} \binom{k-2}{j}\binom{n-1-j}{k} \label{A}\\
& & + (\ell-1)^2\sum_j  (\ell - 2)^{k-j-3} \binom{k-3}{j}\binom{n-2-j}{k-1} \label{B}\\
& & + (\ell-2)(\ell-1)^2 \sum_j  (\ell - 2)^{k-j-3} \binom{k-3}{j}\binom{n-1-j}{k-1}. \label{C}
\end{eqnarray}
In \eqref{A}, \eqref{B} and \eqref{C} we have extended the summations to all integers, noting that in each case the summands are $0$ outside the ranges given by \eqref{fib-rec}.

Bringing the $(\ell-2)$ from outside to inside the summation in \eqref{C}, shifting index from $j$ to $j-1$ in \eqref{B}, and then combining \eqref{B} and \eqref{C}, an application of Pascal's identity shows that \eqref{B} and \eqref{C} sum to
\begin{equation} \label{D}
(\ell-1)^2\sum_j (\ell-2)^{k-j-2}\binom{k-2}{j}\binom{n-1-j}{k-1}.
\end{equation}
Combining \eqref{A} and \eqref{D}, another application of Pascal's identity shows that
$$
p^k_{n,\ell} = (\ell-1)^2\sum_j  (\ell - 2)^{k-j-2} \binom{k-2}{j}\binom{n-j}{k} = (\ell-1)^2\sum_{j=0}^{k-2}  (\ell - 2)^{k-j-2} \binom{k-2}{j}\binom{n-j}{k},
$$
completing the inductive verification of \eqref{pknl-summation-formula}.

\subsection{Combinatorial proof of the formula for  $p_{n,\ell}^k$} \label{subsec-comb-path}

In Section \ref{subsec-induction-path} we established \eqref{pknl-summation-formula} inductively.
Here we give a combinatorial proof of this formula. Let ${\mathcal P}_{n,\ell}^k$ be the set of independent sets of size $k$ in the hypergraph $P_{n,\ell}$. We will partition ${\mathcal P}_{n,\ell}^k$ into $k-1$ blocks, ${\mathcal P}_{n,\ell}^{k,0}, {\mathcal P}_{n,\ell}^{k,1}, \ldots, {\mathcal P}_{n,\ell}^{k,k-2}$, with the property that for each $j$
\begin{equation} \label{block_j_count}
\left|{\mathcal P}_{n,\ell}^{k,j}\right| = (\ell-1)^2(\ell-2)^{k-2-j}\binom{k-2}{j}\binom{n-j}{k}.
\end{equation}
The partitioning will be done by identifying a statistic associated with an independent set (ranging from $0$ to $k-2$); ${\mathcal P}_{n,\ell}^{k,j}$ will be those independent sets with statistic $j$. We will then establish \eqref{block_j_count} by a direct count.

\begin{figure}[htp]
    \centering
    \includegraphics[width=13cm]{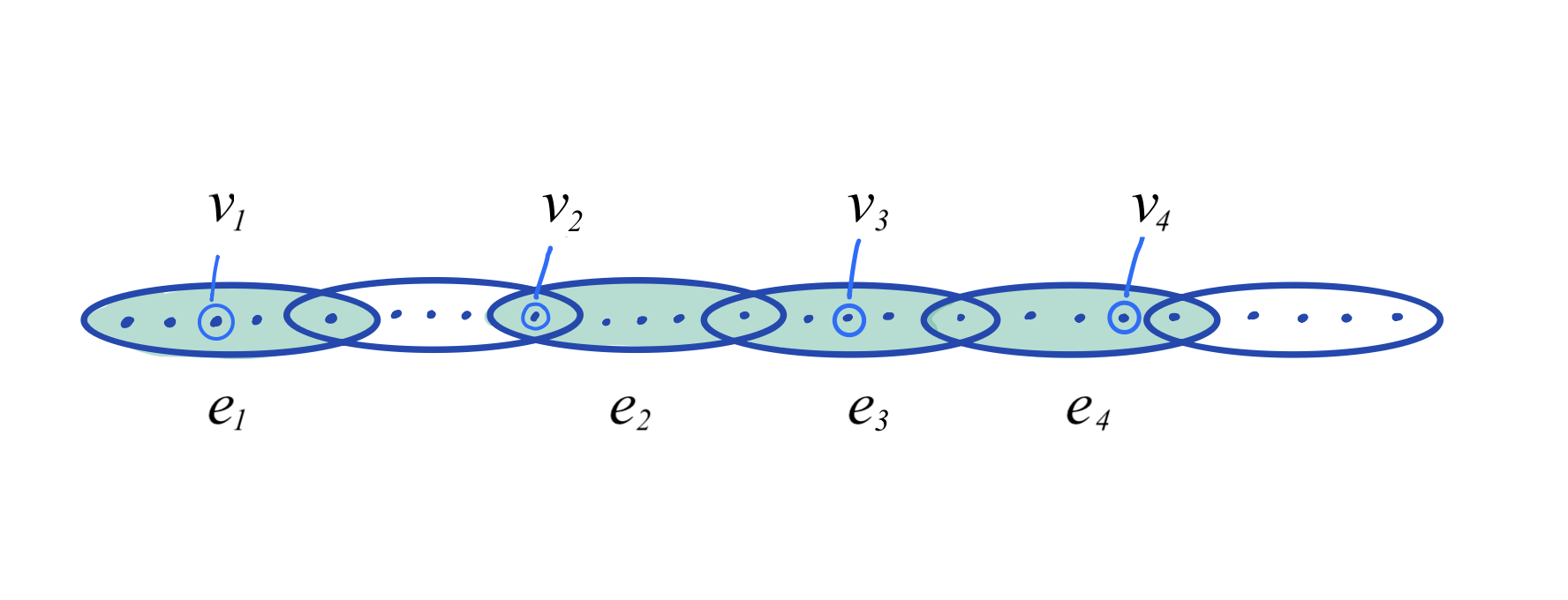}
    \caption{An independent set of size $4$ in the hyperpath $P_{6,5}$. The marked edges ($e_1, e_2, e_3$ and $e_4$) are shaded in teal. The vertex $v_2$ is forced and the vertex $v_3$ is free.}
    \label{fig:exampleskeleton}
\end{figure}

Let $I=\{v_1, \ldots, v_k\}$ be an independent set in $P_{n,\ell}$ of size $k$, with $v_1 < v_2 < \cdots < v_k$ in the natural linear order on the vertices of $P_{n,\ell}$. To $v_i$ we associate a {\em marked} edge $e_i$ of $P_{n,\ell}$ via the following rule: 
\begin{itemize}
\item if there is a unique edge that contains $v_i$, then that unique edge is the marked edge; and
\item if $v_i$ is in two edges of $P_{n,\ell}$, then the later (rightmost) edge of the two is the marked edge. 
\end{itemize}
Note that the number of marked edges is $|I|$. We refer to the collection of marked edges as the {\em skeleton} of $I$. (See Figure \ref{fig:exampleskeleton}.)

For $i=2, \ldots, k-1$, say that $v_i$ is \emph{forced} if $v_i$ is the leftmost vertex of $e_i$ (so, the vertex that $e_i$ has in common with the edge of $P_{n,\ell}$ that comes immediately to its left), and \emph{free} otherwise. (Again see Figure \ref{fig:exampleskeleton}.) Note that for $i=1$ and $i=k$ we do not define forced or free. 

We are now ready to define the statistic $j$:
\begin{quote}
For $j=0, \ldots, k-2$, let ${\mathcal P}_{n,\ell}^{k,j}$ be the set of independent sets in ${\mathcal P}_{n,\ell}^k$ that have exactly $j$ forced vertices.
\end{quote}

\begin{claim} \label{clm-block_j_count}
For each $n \geq 2, k \geq 2, \ell \geq 3$ and $0 \leq j \leq k-2$, \eqref{block_j_count} holds.
\end{claim}

This claim furnishes a combinatorial proof of \eqref{pknl-summation-formula}. 

\begin{proof} (Claim \ref{clm-block_j_count}.)
We will make extensive use of the following observation, which follows from a simple stars-and-bars argument: 
\begin{equation} \label{eq-stars-bars}
\begin{array}{c}
\mbox{The number of solutions to the equation $x_1+x_2+\cdots+x_t = m$, with all $x_i$'s}\\ 
\mbox{non-negative integers (and $m$ a non-negative integer), in which exactly $s$ of}\\
\mbox{the $x_i$'s are constrained to be at least $1$, exactly $r$ of them are constrained }\\
\mbox{to be equal to $0$, and the rest are constrained to be at least $0$, is}\\
\mbox{$\binom{t-r-1+m-s}{t-r-1}$.}
\end{array}
\end{equation}
The relevance of this observation stems from the fact that the $k$ edges in the skeleton give rise to a collection of $k+1$ intervals of consecutive unmarked edges among the $n-k$ unmarked edges. Specifically, the first interval is the collection of unmarked edges to the left of the first marked edge (this collection may be empty). The second interval is the collection of unmarked edges between the first marked edge and the second marked edge, not including the marked edges, and so on, up to the $(k+1)$st interval, which is the collection of unmarked edges to the right of the last marked edge. Letting $a_i$ be the number of unmarked edges in the $i$th interval (so $a_i$ is the number unmarked edges between $e_{i-1}$ and $e_i$ for $i=2, \ldots, k$) we get that each skeleton gives rise to an equation $a_1+a_2+\cdots +a_{k+1} = n-k$, with all the $a_i$'s integers greater than or equal to $0$. 

We further partition ${\mathcal P}_{n,\ell}^{k,j}$ into four blocks:
\begin{description}
\item[${\mathcal P}_{n,\ell}^{k,j}(0,\geq 1)$] is the set of independent sets in ${\mathcal P}_{n,\ell}^{k,j}$ in which $a_k=0$ and $a_{k+1} \geq 1$.
\item[${\mathcal P}_{n,\ell}^{k,j}(0,0)$] is the set of independent sets in ${\mathcal P}_{n,\ell}^{k,j}$ in which $a_k=a_{k+1}=0$.
\item[${\mathcal P}_{n,\ell}^{k,j}(\geq 1,\geq 1)$] is the set of independent sets in ${\mathcal P}_{n,\ell}^{k,j}$ in which $a_k, a_{k+1} \geq 1$.
\item[${\mathcal P}_{n,\ell}^{k,j}(\geq 1,0)$] is the set of independent sets in ${\mathcal P}_{n,\ell}^{k,j}$ in which $a_k \geq 1$ and $a_{k+1}=0$.
\end{description}

We start by enumerating ${\mathcal P}_{n,\ell}^{k,j}(0,\geq 1)$. As described earlier, the skeleton induces a composition $a_1+a_2+\cdots +a_{k+1} = n-k$ with all the $a_i$'s integers greater than or equal to $0$. We have the specific restriction $a_k=0$. 

Additionally, we have $j$ forced vertices. When we force a vertex (say the first vertex of $e_i$), the edge to the left of $e_i$ cannot be marked. Thus the interval of unmarked edges to the left of an edge that has a forced vertex must be non-empty, which means that $j$ of the $a_i$'s must be at least $1$. But we also in this case have unmarked edges following the last marked edge, leading to one more $a_i$ that must be at least $1$ (specifically $a_{k+1}$). Notice the rest of the intervals have no constraints on their lengths\,---\,they may be empty or non-empty.

We may thus apply \eqref{eq-stars-bars} with $m=n-k$, $t=k+1$, $s=j+1$ and $r=1$ to conclude that the number of skeletons in this case is given by $\binom{n-j-2}{k-1}$. 

We must further decide where among the $k-2$ middle edges of the skeleton the $j$ edges that are forced occur, or, equivalently, where the $k-2-j$ free edges occur. This is accounted for by a factor $\binom{k-2}{k-2-j}$.

Inside each of the $k-2-j$ free edges in the skeleton we must choose which vertex from the edge is actually in the independent set. That vertex must be one of the $\ell-2$ vertices in the center of the edge (since the edge is free, we cannot choose the first vertex, and we also cannot choose the last overlapping vertex, as that would be the marked vertex in the next edge). This leads to a factor of $(\ell - 2)^{k-2-j}$.

Considering the first edge $e_1$, we must choose one of the vertices in that edge to be in the independent set. We cannot choose the final vertex of the edge; that would cause the edge immediately to the right of $e_1$ to be marked, and $e_1$ not to be marked. So we have $\ell-1$ options.

The analysis of the previous three paragraphs will be common to all four cases, so it is convenient to set
\begin{equation} \label{A-def}
A:= (\ell-1)(\ell - 2)^{k-2-j}\binom{k-2}{k-2-j}.
\end{equation}

Finally looking at edge $e_k$, the last marked edge: we cannot choose the first vertex of $e_k$ to be in the independent set, or else the previous edge (immediately to the left of $e_k$) could not have been marked. Since the hypergraph has edges after $e_k$, we also cannot choose the last vertex of $e_k$ to be in the independent set (if we did, that vertex would have caused the edge immediately to the right of $e_k$ to have been the marked edge associated with that vertex). Thus we are left with $(\ell-2)$ options to choose from for the vertex of the independent set that is in $e_k$.

It follows that
\begin{equation} \label{block1-01}
\left|{\mathcal P}_{n,\ell}^{k,j}(0,\geq 1)\right| = A(\ell-2)\binom{n-j-2}{k-1}.
\end{equation}

We now move on to enumerating ${\mathcal P}_{n,\ell}^{k,j}(0,0)$. Again the skeleton induces a composition $a_1+a_2+\cdots +a_{k+1} = n-k$ with all the $a_i$'s integers greater than or equal to $0$, and with the specific restriction $a_k=0$. 

As before each of the $j$ forced vertices gives rise to a restriction $a_i \geq 1$, and in this case there is no further such restriction. We may thus apply \eqref{eq-stars-bars} with $m=n-k$, $t=k+1$, $s=j$ and $r=1$ to conclude that the number of skeletons in this case is given by $\binom{n-j-2}{k-2}$. 

Deciding where among the $k-2$ middle edges of the skeleton the $j$ edges that are forced occur, choosing a vertex from inside each of the $k-2-j$ free edges to be in the independent set, and choosing a vertex from $e_1$ to be in the independent set, together give a contribution of $(\ell-1)(\ell - 2)^{k-2-j}\binom{k-2}{k-2-j}$, as before. 

To finish this case, again we look at the last marked edge $e_k$. Again we cannot choose the first vertex to be a part of the independent set, or else the previous edge could not be marked (contradicting $a_k=0$). Since the last vertex of $e_k$ is not the first vertex of a later edge, we can choose the final vertex of $e_k$ to be the vertex of the independent set from $e_k$. Thus in this case we have $(\ell-1)$ options to choose from for the vertex of the independent set that is in $e_k$.

It follows that
\begin{equation} \label{block2-00}
\left|{\mathcal P}_{n,\ell}^{k,j}(0,0)\right| = A(\ell-1)\binom{n-j-2}{k-2}.
\end{equation}

To enumerate ${\mathcal P}_{n,\ell}^{k,j}(\geq 1,\geq 1)$, note that here we have $a_k \geq 1$ instead of $a_k=0$, but otherwise our considerations of $r$ and $s$ are the same as in the case of ${\mathcal P}_{n,\ell}^{k,j}(0,\geq 1)$. Thus $s = j+2$ and $r=0$ and this leads to a count of $\binom{n-j-2}{k}$ for the number of skeletons.

Deciding where among the $k-2$ middle edges of the skeleton the $j$ edges that are forced occur, choosing a vertex from inside each of the $k-2-j$ free edges to be in the independent set, and choosing a vertex from $e_1$ to be in the independent set, together give a contribution of $(\ell-1)(\ell - 2)^{k-2-j}\binom{k-2}{k-2-j}$, as before. 

Now looking at $e_k$, we find that, unlike the first two cases, we can include the first vertex as the vertex from $e_k$ that causes it to be marked, since the edge immediately to the left of $e_k$ is unmarked. And since there are edges in the hyperpath after (to the right of) $e_k$ in this case, similar to our consideration of ${\mathcal P}_{n,\ell}^{k,j}(0,\geq 1)$ we cannot chose the final vertex of $e_k$ to be the vertex of $e_k$ that is in the independent set. Thus we have $(\ell-1)$ options for the vertex of the independent set coming from $e_k$.

It follows that
\begin{equation} \label{block3-11}
\left|{\mathcal P}_{n,\ell}^{k,j}(\geq 1,\geq 1)\right| = A(\ell-1)\binom{n-j-2}{k}.
\end{equation}

Finally, we consider ${\mathcal P}_{n,\ell}^{k,j}(\geq 1,0)$. Similar to arguments in the previous cases, we have here $r=1$ and $s=j+1$, so we have $\binom{n-j-2}{k-1}$ as the count of the number of skeletons. 

As in all previous cases, deciding where among the $k-2$ middle edges of the skeleton the $j$ edges that are forced occur, choosing a vertex from inside each of the $k-2-j$ free edges to be in the independent set, and choosing a vertex from $e_1$ to be in the independent set, together give a contribution of $(\ell-1)(\ell - 2)^{k-2-j}\binom{k-2}{k-2-j}$. 

Finally we look again to edge $e_k$. Similarly to our consideration of ${\mathcal P}_{n,\ell}^{k,j}(\geq 1,\geq 1)$ we can choose the first vertex of $e_k$ to be the vertex from $e_k$ that is in the independent set, and similarly to ${\mathcal P}_{n,\ell}^{k,j}(0,0)$ we can also choose the final vertex. Thus we have $\ell$ options in this case for the final vertex of the independent set.

It follows that
\begin{equation} \label{block4-10}
\left|{\mathcal P}_{n,\ell}^{k,j}(\geq 1,\geq 1)\right| = A \ell \binom{n-j-2}{k-1}.
\end{equation}

Combining \eqref{block1-01}, \eqref{block2-00}, \eqref{block3-11} and \eqref{block4-10} we have   
\begin{eqnarray}
\frac{\left|{\mathcal P}_{n,\ell}^{k,j}\right|}{A} & = & (\ell-2)\binom{n-j-2}{k-1} + (\ell-1)\binom{n-j-2}{k-2} + \nonumber\\
& & ~~~~~~~~~~~~~~~~~~~~~~~~~~~~~~~~~~~~~~~~~~(\ell-1)\binom{n-j-2}{k} + \ell\binom{n-j-2}{k-1} \nonumber \\
& = & (\ell-1)\left(\binom{n-j-2}{k-1} + \binom{n-j-2}{k-2} + \binom{n-j-2}{k} + \binom{n-j-2}{k-1}\right) \label{exp1} \\
& = & (\ell-1)\left(\binom{n-j-1}{k-1} + \binom{n-j-1}{k}\right) \label{exp2} \\
& = & (\ell-1)\binom{n-j}{k} \label{exp3}.
\end{eqnarray}

Here \eqref{exp1} comes from a simple rearrangement of summands. Specifically  we take one of the $\binom{n-j-2}{k-1}$'s from the $\ell\binom{n-j-2}{k-1}$ term and move it over to the $(\ell-2)\binom{n-j-2}{k-1}$ term. We also use Pascal's identity twice in \eqref{exp2} on adjacent binomial coefficients and again in \eqref{exp3}. Recalling the definition of $A$ from \eqref{A-def} we obtain
$$
\left|{\mathcal P}_{n,\ell}^{k,j}\right| = (\ell-1)(\ell-2)^{k-2-j}\binom{n-j}{k}\binom{k-2}{k-2-j},
$$
which is \eqref{block_j_count}.
\end{proof}

\subsection{Proofs of \eqref{fib-rec} and \eqref{eq-Px-pure}} \label{subsec-fib-rec}

In Section \ref{subsec-induction-path} we derived \eqref{fib-rec} by algebraic manipulation. Such a simple formula should admit a combinatorial explanation (and indeed, as we have observed, there is a very simple combinatorial explanation in the case $\ell=2$). 

Here we present a bijection, due to Ferenc Bencs (see Note \ref{note-comb}), that directly explains \eqref{fib-rec}. Recall that $P_{n,\ell}$ has consecutive edges $e_1, \ldots, e_n$. Let 
\begin{itemize} \item $v$ be the unique vertex in common to $e_{n-1}$ and $e_n$,
\item $v_1, \ldots, v_{\ell-2}$ be the vertices of $e_{n-1}$ that are only in $e_{n-1}$, 
\item $w$ be a vertex chosen arbitrarily from $e_n\setminus \{v\}$, and
\item $w_1, \ldots, w_{\ell-2}$ be the vertices of $e_n\setminus \{v,w\}$.
\end{itemize}

Let ${\mathcal I}^k_{n,\ell}$ be the set of independent sets in $P_{n,\ell}$ of size $k$. We partition ${\mathcal I}_{n,\ell}$ into $2\ell-2$ blocks, falling into four categories (the first two consisting of one block each, the third and fourth consisting of $\ell-2$ blocks each):
\begin{enumerate}
\item ${\mathcal A}^k_{n,\ell}$ consists of those independent sets in ${\mathcal I}^k_{n,\ell}$ that include $w$ (and so necessarily do not include any of $w_1, \ldots, w_{\ell-2}$ or $v$) and include none of $v_1, \ldots, v_{\ell-2}$. There is a simple bijection from ${\mathcal A}^k_{n,\ell}$ to ${\mathcal I}^{k-1}_{n-2,\ell}$, obtained by deleting $w$ from an independent set in ${\mathcal A}^k_{n,\ell}$. It follows that $|{\mathcal A}^k_{n,\ell}|=p^{k-1}_{n-2,\ell}$.
\item ${\mathcal B}^k_{n,\ell}$ consists of those independent sets in ${\mathcal I}^k_{n,\ell}$ that include none of $w$ or $w_1, \ldots, w_{\ell-2}$. The identity map is a bijection from ${\mathcal B}^k_{n,\ell}$ to ${\mathcal I}^k_{n-1,\ell}$, so $|{\mathcal B}^k_{n,\ell}|=p^k_{n-1,\ell}$.
\item For each $i=1, \ldots, \ell-2$, ${\mathcal C}^k_{n,\ell}(i)$ consists of those independent sets in ${\mathcal I}^k_{n,\ell}$ that include $w_i$ (and so necessarily do not include any of $v, w$ or $w_j$, $j \neq i$. The identity map is a bijection from ${\mathcal C}^k_{n,\ell}(i)$ to 
$$
\{I \in {\mathcal I}^{k-1}_{n-1,\ell}: v \not \in I\}.
$$
\item For each $i=1, \ldots, \ell-2$, ${\mathcal D}^k_{n,\ell}(i)$ consists of those independent sets in ${\mathcal I}^k_{n,\ell}$ that include both $w$ and $v_i$ (and so necessarily do not include any of $v$, $w_1, \ldots, w_{\ell-2}$, or $v_j$, $j \neq i$. There is a simple bijection from ${\mathcal D}^k_{n,\ell}(i)$ to 
$$
\{I \in {\mathcal I}^{k-1}_{n-1,\ell}: v \in I\},
$$
obtained by deleting $w$ and $v_i$ from an independent set in ${\mathcal D}^k_{n,\ell}(i)$, and adding $v$.
\end{enumerate}
Combining the observations in points 3 and 4 above we obtain
$$
\left|\cup_{i=1}^{\ell-2} \left({\mathcal C}^k_{n,\ell}(i) \cup {\mathcal D}^k_{n,\ell}(i)\right)\right| = (\ell-2)p^{k-1}_{n-1,\ell}.
$$
Combining this with the observations in the points 1 and 2 above we obtain
$$
p^k_{n,\ell} = p^k_{n-1,\ell}+p^{k-1}_{n-2,\ell}+(\ell-2)p^{k-1}_{n-1,\ell}.
$$

\subsection{The independent set sequence of hypercombs and related hypergraphs} \label{subsec-comb-proofs}

We begin by describing the more general setting in which we will work, alluded to just after Note \ref{note-hcomb}. Let $G$ be a hypergraph (not necessarily uniform, not necessarily a hypertree) with distinguished vertex $v$. For $\ell \geq 3$, let a sequence $(H_{n,\ell})_{n \geq 0}$ of hypergraphs be defined as follows:
\begin{itemize}
\item $H_{0,\ell}$ consists of a single copy of $G$.
\item $H_{1,\ell}$ consists of an edge of size $\ell$, together with two disjoint copies of $G$; in the first of these copies, the distinguished vertex $v$ is identified with some (arbitrary) vertex of the initial edge, and in the second, it is identified with some other (arbitrary) vertex of the initial edge. 
\item For $n \geq 2$, $H_{n,\ell}$ starts with a copy of $P_{n,\ell}$. Let $w_0$ be a vertex in $e_1 \setminus e_2$. For $i=1,\ldots, n-1$ let $w_i$ be the vertex in common to $e_i$ and $e_{i+1}$. Finally let $w_n$ be a vertex in $e_n \setminus e_{n-1}$. Augment $P_{n,\ell}$ by adding $n+1$ disjoint copies of $G$, say $G^0, \ldots, G^n$, with the distinguished vertex in $G^i$ identified with $w_i$ for $i=0, \ldots, n$.  
\end{itemize}
Notice that when $G$ is a single edge of order $\ell$, $H_{n,\ell}$ coincides with $C_{n,\ell}$ for $n \geq 1$. Note also that in this case the independence polynomial of $H_{0,\ell}$ is $1+\ell x$, exactly what $C_{0,\ell}(x)$ was declared to be in Theorem \ref{thm-comb}. 

Let also a sequence $(D_{n,\ell})_{n \geq 0}$ of auxiliary hypergraphs be defined as follows:
\begin{itemize}
\item $D_{n,\ell}$ is obtained from $H_{n,\ell}$ by removing the first (left-most) copy of $G$ from $H_{n,\ell}$ (including the vertex at which the copy of $G$ meets the underlying copy of $P_{n,\ell}$).
\end{itemize}
In particular this means that $D_{0,\ell}$ is empty (so $D_{0,\ell}(x)=1$). An illustration of $D_{n,\ell}$ when $G$ is a single edge of order $\ell$ is given in Figure \ref{fig:CandD}.

Let $H_{n,\ell}(x)$ and $D_{n,\ell}(x)$ denote the independence polynomials of $H_{n,\ell}$ and $D_{n,\ell}$ respectively. Let also 
\begin{itemize}
\item $G_0(x)$ be the independence polynomial of $G$,
\item $G_1(x)$ be the independence polynomial of the hypergraph obtained from $G$ by deleting the vertex $v$ (which reduces by $1$ the order of every edge in $G$ that includes $v$), and
\item $G_2(x)$ be the independence polynomial of the hypergraph obtained from $G$ by deleting the vertex $v$ and all edges in $G$ that include $v$.
\end{itemize}

We now observe that there are coupled recurrence relations for $H_{n,\ell}(x)$ and $D_{n,\ell}(x)$ valid for $n \geq 1$. For $H_{n,\ell}(x)$ we have
\begin{eqnarray}
H_{n,\ell}(x) & = & G_1(x)D_{n,\ell}(x)+xG_1(x)G_2(x)D_{n-1,\ell}(x) \nonumber \\
& = & G_1(x)\left(D_{n,\ell}(x)+xG_2(x)D_{n-1,\ell}(x)\right) \label{rec-H}
\end{eqnarray}
(consider first not occupying the first vertex along the spine of $H_{n,\ell}(x)$ that has a copy of $G$ dropped from it, and then occupying it), and for $D_{n,\ell}(x)$ we have
\begin{equation} \label{rec-D}
D_{n,\ell}(x)=H_{n-1,\ell}(x)+x(\ell-2)G_1(x)D_{n-1,\ell}(x)
\end{equation}
(consider first not occupying any of the initial $\ell-2$ vertices along the spine that do not have a copy of $G$ dropped from them, then occupying one of them). Applying \eqref{rec-D} twice (once with index $n$ and once with index $n-1$) and rearranging terms we get
\begin{eqnarray*}
D_{n,\ell}(x)+xG_2(x)D_{n-1,\ell}(x) & = & H_{n-1,\ell}(x)+x(\ell-2)G_1(x)D_{n-1,\ell}(x) \\
& & + xG_2(x)H_{n-2,\ell}(x)+x^2(\ell-2)G_1(x)G_2(x)D_{n-2,\ell}(x) \\
& = & (1+(\ell-2)x)H_{n-1,\ell}(x)+xG_2(x)H_{n-2,\ell}(x).
\end{eqnarray*}
In the second equality above we use \eqref{rec-H} to combine the second and fourth terms on the right-hand side. It follows (using \eqref{rec-H} again) that
\begin{equation} \label{H-rec}
H_{n,\ell}(x) = G_1(x)\left((1+(\ell-2)x)H_{n-1,\ell}(x)+xG_2(x)H_{n-2,\ell}(x)\right).
\end{equation}
Noting that when $G$ is a single edge of order $\ell$ we have $G_1(x)=1+(\ell-1)x$ and $G_2(x)=1$ (since $G_2$ is empty in this case), we obtain \eqref{eq-comb-rec} as a special case of \eqref{H-rec}.  

The initial conditions for the recurrence \eqref{H-rec} are 
$$
H_{0,\ell}(x) = G_0(x)~~~\mbox{and}~~~H_{1,\ell}(x) = G_0(x)^2 -x^2G_2(x)^2 + (\ell-2)xG_1(x)^2.
$$ 
For the second of these initial conditions, consider first those independent sets in $H_{1,\ell}$ that only use vertices from the two copies of $G$ in $H_{1,\ell}$. The factor $G_0(x)^2$ comes from choosing an arbitrary independent set from each copy of $G$. The factor of $-x^2G_2(x)^2$ comes from the fact that we have over-counted by including all pairs that have the property that each independent set in the pair includes a vertex from the spine. Then consider those independent sets in $H_{1,\ell}$ that include one vertex from among the $\ell-2$ vertices that are not part of a copy of $G$. These sets contribute $(\ell-2)xG_1(x)^2$ to the independence polynomial. When $G$ is a single edge of order $\ell$ these formulae are easily seen to give the initial conditions posited in Theorem \ref{thm-comb} for the hypercomb, completing the proof of that theorem. 

We now move on to Corollary \ref{cor-comb}, the partial unimodality of the independent set sequence of $C_{n,\ell}$. Note that the graph obtained from $C_{n,\ell}$ by replacing each edge with a clique is not claw-free for $n \geq 2$, so we cannot apply the results of Hamidoune and of Chudnovsky and Seymour that we used for the linear hyperpath. 

We work in the general setting of a sequence $(z_n(x))_{n \geq 0}$ of polynomials with non-negative coefficients defined by the recurrence
$$
z_n(x) = a(x)z_{n-1}(x) + b(x)z_{n-2}(x)
$$
for $n \geq 2$, where $a(x), b(x)$ are non-negative and non-zero polynomials. Our goal is to establish easy-to-check conditions on $a(x)$, $b(x)$, $z_0(x)$ and $z_1(x)$ that imply unimodality of all but a vanishing proportion of the coefficient sequence of $z_n(x)$. We briefly sketch the approach here.
\begin{itemize}
\item First we explicitly solve the recurrence, to express $z_n(x)$ as the sum of two polynomials, $W_n(x)$ and $Y_n(x)$, each of which can be explicitly factored into a product of low-degree polynomials. 
\item Assuming that each of these factors is log-concave we conclude that $W_n(x)$ and $Y_n(x)$ are log-concave and so unimodal. 
\item Viewing $W_n$ and $Y_n$ as probability generating functions of discrete distributions we explicitly compute the expectations of $W_n$ and $Y_n$ and find that they are (asymptotically) the same.
\item We then appeal to a result of Bottomley, that shows that the expectation and the mode of a unimodal polynomial are (quantifiably) close to each other. We conclude that the modes of $W_n$ and $Y_n$ are close to each other.
\item  Since the sum of two unimodal polynomials is increasing at least up to the smaller of the two modes, and is decreasing at least from the larger of the two modes on, we conclude that $z_n(x)=W_n(x)+Z_n(x)$ is mostly unimodal (specifically, unimodal except possibly between the modes of $W_n(x)$ and $Y_n(x)$).  
\end{itemize}
This scheme leads to a general statement (Theorem \ref{thm-almost-unimodality} below) to the effect that if a particular collection of low-degree polynomials is log-concave, then $z_n(x)$ is mostly unimodal for all $n$. We then apply this general result to the linear uniform hypercomb.  

Let us now begin executing the above-described scheme. In the initial part of what follows we draw on the presentation given in \cite{WZ2011}. First note that using standard recurrence relation techniques we have that for $n \geq 2$
\begin{eqnarray}
z_n(x) & = & \frac{(z_1(x)-z_0(x)\mu)\lambda(x)^n + (z_0(x)\lambda(x)-z_1(x))\mu(x)^n}{\lambda(x)-\mu(x)} \nonumber \\
& = & z_1(x)\left(\frac{\lambda(x)^n-\mu(x)^n}{\lambda(x)-\mu(x)}\right) -z_0(x)\lambda(x)\mu(x)\left(\frac{\lambda(x)^{n-1}-\mu(x)^{n-1}}{\lambda(x)-\mu(x)}\right) \label{gen-rec-sol}
\end{eqnarray}
where
$$
\lambda(x) = \frac{a(x)+\sqrt{a^2(x)+4b(x)}}{2}~~~\mbox{and}~~~\mu(x) = \frac{a(x)-\sqrt{a^2(x)+4b(x)}}{2}
$$
or equivalently
$$
\lambda(x)+\mu(x) = a(x)~~~\mbox{and}~~~\lambda(x)\mu(x) = -b(x).
$$
(See e.g. \cite[Chapter 7]{B2010}.) Note that since $a(x)$, $b(x)$ are non-negative and non-zero we have $a^2(x)+4b(x)>0$ for $x>0$ and so $\lambda(x) \neq \mu(x)$.

Next (see e.g. \cite[Chapter V.19]{BC1955}), we have the following factorizations for arbitrary $p, q$ with $p \neq q$:
$$
\frac{p^n-q^n}{p-q} = \left\{
\begin{array}{cc}
\prod_{s=1}^{(n-1)/2} \left((p+q)^2 -4pq\cos^2 \left(\frac{s\pi}{n}\right) \right) & \mbox{if $n \geq 1$ is odd,}\\
(p+q)\prod_{s=1}^{(n-2)/2} \left((p+q)^2 -4pq\cos^2 \left(\frac{s\pi}{n}\right) \right) & \mbox{if $n \geq 2$ is even}
\end{array}
\right.
$$
(here interpreting the empty product to be $1$, as usual).
Combining this with \eqref{gen-rec-sol} we get that
$$
z_n(x) = W_n(x)+Y_n(x)
$$
for $n \geq 2$, where
\begin{eqnarray} \label{W-def}
W_n(x) & = & 
\left\{
\begin{array}{rl}
z_1(x)a(x)\prod_{s=1}^{(n-2)/2} \left(a^2(x) + 4b(x)\cos^2\left(\frac{s\pi}{n}\right)\right) & \mbox{if $n$ is even,} \\
z_1(x)\prod_{s=1}^{(n-1)/2} \left(a^2(x) + 4b(x)\cos^2\left(\frac{s\pi}{n}\right)\right) & \mbox{if $n$ is odd}
\end{array}
\right.
\end{eqnarray}
and 
\begin{eqnarray} \label{Y-def}
Y_n(x) & = & 
\left\{
\begin{array}{rl}
z_0(x)b(x)\prod_{s=1}^{(n-2)/2} \left(a^2(x) + 4b(x)\cos^2\left(\frac{s\pi}{n-1}\right)\right) & \mbox{if $n$ is even,} \\
z_0(x)a(x)b(x)\prod_{s=1}^{(n-3)/2} \left(a^2(x) + 4b(x)\cos^2\left(\frac{s\pi}{n-1}\right)\right) & \mbox{if $n$ is odd.}
\end{array}
\right.
\end{eqnarray}
Let us now assume that each of the following polynomials:
$$
z_0(x), z_1(x), a(x), b(x)~\mbox{and}~ a^2(x)+4b(x)y~(y \in [0,1])
$$
have log-concave coefficient sequences. Using the standard fact that the log-concavity of coefficient sequences of polynomials is preserved under multiplication (see e.g. \cite[Theorem 9.9]{B2016}), it follows that each of $W_n(x)$ and $Y_n(x)$ have log-concave and so unimodal coefficient sequences.  

Let $mw_n$ be a mode of $W_n(x)$. That is, let $mw_n$ be an integer satisfying that the coefficient sequence of $W_n(x)$ is weakly increasing up to (and including) the coefficient of $x^{mw_n}$ and is weakly decreasing starting from the coefficient of $x^{mw_n}$. Let $my_n$ be a mode of $Y_n(x)$. Evidently the coefficient sequence of $z_n(x)$ is weakly increasing at least up to the coefficient of $x^{\min\{mw_n,my_n\}}$, and it is weakly decreasing at least from the coefficient of $x^{\max\{mw_n,my_n\}}$ on. So the smaller we can make $|mw_n-my_n|$, the closer we get to concluding unimodality of the coefficient sequence of $z_n(x)$. 

Let us focus first on estimating $mw_n$ in the case where $n$ is even. Because $W_n(x)$ is a non-negative and non-zero polynomial we may view $W_n(x)/W_n(1)$ as the probability generating function of a discrete probability distribution $X_n$, supported on a finite subset of the natural numbers. In light of the factorization we have found for $W_n(x)$ we see that 
$$
X_n = X_n^0 + \sum_{s=1}^{(n-2)/2} X_n^s
$$
where $X_n^0$ has probability generating function $z_1(x)a(x)/(z_1(1)a(1))$, $X_n^s$ has probability generating function $(a^2(x) + 4b(x)\cos^2(s\pi/n))/((a^2(1) + 4b(1)\cos^2(s\pi/n)))$ for each $s$, and the $X_n^i$'s are independent. 

We now use that if a distribution $X$ has probability generating function $P_X(x)$ then its expectation is given by $E(X) = P'_X(x)|_{x=1}$, where $'$ indicates derivative with respect to $x$. So 
$$
E(X_n^0)=\frac{(z_1(x)a(x))'|_{x=1}}{z_1(1)a(1)} := c_0
$$ 
and for $s \geq 1$
$$
E(X_n^s) = \frac{c_1 + c_2\cos^2(s\pi/n)}{c_3+c_4\cos^2(s\pi/n)} $$
where the $c_i$'s are constants that can be explicitly calculated from $a(x)$, $b(x)$ and (in the case of $c_0$) $z_1(x)$, with $c_3, c_4 \neq 0$. It follows that
\begin{equation} \label{expectation}
E(X_n) = c_0 + \sum_{s=1}^{(n-2)/2} \frac{c_1 + c_2\cos^2(s\pi/n)}{c_3+c_4\cos^2(s\pi/n)}.
\end{equation}
We estimate the sum in \eqref{expectation} using an integral. Setting $m=(n-2)/2$ we have
\begin{equation} \label{sum-translation}
\sum_{s=1}^{(n-2)/2} \frac{c_1 + c_2\cos^2(s\pi/n)}{c_3+c_4\cos^2(s\pi/n)} = m\sum_{s=1}^m \frac{1}{m}\left(\frac{c_1 + c_2\cos^2((2s/(2m+2))(\pi/2))}{c_3+c_4\cos^2((2s/(2m+2))(\pi/2))}\right). 
\end{equation}
Observing that 
\begin{equation} \label{Rs-check}
\frac{s-1}{m} \leq \frac{2s}{2m+2} \leq \frac{s}{m}
\end{equation}
(the latter is always true, the former is true for $s \leq m+1$) we see that the summation on the right-hand side of \eqref{sum-translation} is a Riemann sum estimate for $\int_0^1 f(t)~dt$ where
$$
f(t) = \frac{c_1+c_2\cos^2(t\pi/2)}{c_3+c_4\cos^2(t\pi/2)}.
$$
Since $f'(t)$ is bounded above and below on the interval $[0,1]$ by a constant depending on $c_1$, $c_2$, $c_3$ and $c_4$, using standard Riemann sum error estimates we get
$$
\sum_{s=1}^m \frac{1}{m}\left(\frac{c_1 + c_2\cos^2((2s/(2m+2))(\pi/2))}{c_3+c_4\cos^2((2s/(2m+2))(\pi/2))}\right) = \int_0^1 f(t)~dt \pm \frac{c_5}{m} 
$$
where $c_5$ is a constant depending $c_1$, $c_2$, $c_3$ and $c_4$. (Here we use the notation $A=B \pm C$ as shorthand for $|A-B| \leq C$.) Returning to \eqref{expectation} and using that $m=n/2+O(1)$ we conclude that (for even $n$) there is a constant $K$ (depending only on $a(x)$, $b(x)$) such that
$$
E(X_n) = Kn + O(1),
$$
where the implicit constant in the $O(1)$ depends on $a(x)$, $b(x)$ and $z_1(x)$. Specifically we have
$$
K=\frac{1}{2}\int_0^1 \left(\frac{c_1+c_2\cos^2(t\pi/2)}{c_3+c_4\cos^2(t\pi/2)}\right)dt.
$$

Next we estimate the variance of $X_n$. By independence of the $X_n^i$'s we have 
$$
{\rm Var}(X_n) = {\rm Var}(X_n^0) + \sum_{s=1}^{(n-2)/2} {\rm Var}(X_n^s).
$$ 
We now use that if a distribution $X$ has probability generating function $P_X(x)$ then
$$
{\rm Var}(X) = P''_X(x)|_{x=1} + P'_X(x)|_{x=1} - (P'_X(x)|_{x=1})^2. 
$$
From this it is immediate that ${\rm Var}(X_n^0)$ is a constant (depending on $a(x)$, $b(x)$ and $z_1(x)$), and that for $s \geq 1$ ${\rm Var}(X_n^s)$ can be bounded above and below by constants independent of $s$ (depending only on $a(x)$ and $b(x)$). It follows that $\sigma(X_n)$, the standard deviation of $X_n$, satisfies
$$
\sigma(X_n) \leq O(\sqrt{n})
$$
where the implicit constant in the $O(\sqrt{n})$ depends only on $a(x)$ and $b(x)$.

\medskip

We now appeal to a result of Bottomley \cite{B2004}.
\begin{thm} \label{thm-bot}
Let $X$ be a discrete probability distribution supported on the integers, that is unimodal in the sense that there is an integer $m(X)$ (not necessarily unique) with $P(X=i-1) \leq P(X=i)$ for all $i \leq m(X)$ and $P(X=i) \geq P(X=i+1)$ for all $i \geq m(X)$. If $X$ has finite expectation $E(X)$ and finite standard deviation $\sigma(X)$ then
$$
\left|E(X)-m(X)\right| \leq \sqrt{3} \sigma(X).
$$
\end{thm}
This may be viewed as an extension to arbitrary polynomials of a weak form of Darroch's theorem \cite{D1964}, which recall says that if $X$ is finitely supported on the natural numbers and if furthermore the probability generating polynomial of $X$ has all real roots (a stronger condition than unimodality) then $\left|E(X)-m(X)\right| \leq 1$.

Applying Theorem \ref{thm-bot} to $X_n$ (when $n$ is even) we conclude that
\begin{equation} 
\label{eq-Wmode}
\left|mw_n-Kn\right| \leq O(\sqrt{n}),
\end{equation}
where the implicit constant in the $O(\sqrt{n})$ depends only on $a(x)$ and $b(x)$.

We can do an almost identical analysis for $Y_n$ when $n$ is even. In this case we still set $m=(n-2)/2$ to get the analog of \eqref{sum-translation}. To confirm that the summation on the right-hand side of the analog of \eqref{sum-translation} is still a Riemann sum estimate for $\int_0^1 f(t)~dt$, we need to check an analog of \eqref{Rs-check}, namely
$$
\frac{s-1}{m} \leq \frac{2s}{2m+1} \leq \frac{s}{m}.
$$
We conclude (routine analysis, very similar to that used in the derivation of \eqref{eq-Wmode}, omitted) that 
\begin{equation} 
\label{eq-Ymode}
\left|my_n-Kn\right| \leq O(\sqrt{n})
\end{equation}
for (crucially) the same constant $K$ as in \eqref{eq-Wmode}; and running through the analysis for $n$ odd we find that \eqref{eq-Wmode} and \eqref{eq-Ymode} both hold for odd $n$, again with the same constant $K$. We conclude that the modes of $Y_n$ and $W_n$ differ by at most $O(\sqrt{n})$, and so the coefficient sequence of $z_n(x)$ is unimodal except (possibly) for a $O(\sqrt{n})$ portion around the $[Kn]$th term.   

\medskip

We summarize all of this in a theorem.
\begin{thm} \label{thm-almost-unimodality}
Let $a(x)$ and $b(x)$ be polynomials with non-negative coefficients, each with at least one positive coefficient. Let $z_0(x)$ and $z_1(x)$ be polynomials with non-negative coefficients. For $n \geq 2$ let the polynomial $z_n(x)$ be defined by the recurrence
$$
z_n(x) = a(x)z_{n-1}(x) + b(x)z_{n-2}(x).
$$

Suppose that the coefficient sequences of each of $a(x)$, $b(x)$, $z_0(x)$ and $z_1(x)$ are log-concave, and also that the coefficient sequence of $a^2(x) + 4b(x)y$ is log-concave for all choices of $y \in [0,1]$. Then there are constants $K > 0$ and $L > 0$, both depending only the coefficients of $a(x)$ and $b(x)$, such that for all $n$ the coefficient sequence of $z_n(x)$ is
\begin{itemize}
\item weakly increasing up to the coefficient of $x^{[Kn-L\sqrt{n}]}$ and
\item weakly decreasing from the coefficient of $x^{[Kn+L\sqrt{n}]}$ on.
\end{itemize}
\end{thm}

\medskip

\begin{note}
We can strengthen Theorem \ref{thm-almost-unimodality} very slightly by observing that we do not need $W_n(x)$ and $Y_n(x)$ to be log-concave, just unimodal, and then using the fact that the product of a unimodal polynomial and a log-concave polynomial is unimodal. So we can allow up to one of the factors in each of the four factorizations in \eqref{W-def}, \eqref{Y-def} to be unimodal but not necessarily log-concave, without changing the conclusion of the theorem.   
\end{note}

\medskip

We now apply Theorem \ref{thm-almost-unimodality} to the independent set sequence of the family of linear hypercombs, in order to prove Corollary \ref{cor-comb}. In this case we have 
$$
a(x) = (1+(\ell-1)x)(1+(\ell-2)x),
$$
$$
b(x) = x(1+(\ell-1)x)
$$
and
$$
z_0(x) = C_{0,\ell}(x) = 1 + \ell x,
$$
all three of which are easily seen to have log-concave independent set sequences for all $\ell$ (for $a(x)$ and $b(x)$ we can use that log-concavity is preserved under multiplication, and that linear polynomials are trivially log-concave; for $z_0(x)$ we need only use this second fact). We have
$$
z_1(x) = C_{0,\ell}(x) = 1 + (3\ell-2)x + (2(\ell-1)(\ell-2)+(\ell-1)^2 + 2(\ell-1))x^2 + (\ell-1)^2(\ell-2)x^3. 
$$
The two log-concavity relations that need to be checked here reduce to
$$
6\ell^2 -6\ell +1 \geq 0
$$
and
$$
6\ell^4-22\ell^3+31\ell^2-20\ell+5 \geq 0,
$$
both of which are easily checked to hold for $\ell \geq 3$. Finally we have
$$
a^2(x) + 4b(x)y = (1+(\ell-1)x)(1+(3\ell-5+4y)x + (3\ell^2-10\ell+8)x^2 + (\ell^3-5\ell^2+8\ell-4)x^3).
$$
The first factor trivially has log-concave coefficient sequence. There are two relations that need to be checked to verify that the second factor has log-concave coefficient sequence, and these simplify to
$$
6\ell^2 -(20-24y)\ell + +(17-40y) \geq 0 
$$
and
$$
6\ell^4 -(40+4y)\ell^3 + (99+20y)\ell^2 - (108+32y)\ell + (44+16y) \geq 0. 
$$
We appeal to {\tt Mathematica} to minimize the left-hand side of both inequalities over $\ell\in[3,\infty]$ and $y \in [0,1]$ and see that the minima are both positive, but this could in principle also be done by hand. 

The conditions of Theorem \ref{thm-almost-unimodality} having been verified, we conclude that for each $\ell \geq 3$ there are constants $K_\ell, L_\ell >0$ such that the independent set sequence of the $\ell$-uniform linear hypercomb built on an $n$-edge spine is weakly increasing up to $K_\ell n - L_\ell\sqrt{n}$ and weakly decreasing from $K_\ell n + L_\ell\sqrt{n}$, as asserted in Corollary \ref{cor-comb}.  

\section*{Acknowledgment}

We thank Ferenc Bencs for very helpful comments on an early version of this paper.

\end{document}